\documentclass[12pt]{article}
\usepackage{comment}
\usepackage{url}
\usepackage{amsmath}
\usepackage{amssymb}  
\usepackage{graphicx}
\usepackage{float}
\usepackage{hyperref}
\hypersetup{
    colorlinks=true,
    linkcolor=blue,
    filecolor=magenta,      
    urlcolor=blue,
    pdftitle={Overleaf Example},
    pdfpagemode=FullScreen,
    }

\newcommand{\Z}{{\sf Z}}
\newcommand{\N}{{\mathbb{N}}}
\newcommand{\Q}{{\mathbb{Q}}}
\newcommand{\D}{{\mathbb{D}}}
\renewcommand{\O}{{\mathbb{O}}}

\newcommand{\Ad}{{\mathbb{A_{\hbox{$d$}}}}}
\newcommand{\Af}{{\mathbb{A_{\hbox{$4$}}}}}

\newcommand{\s}[1]{\s_{#1}}

\newcommand{\monus}{\;\raise.5ex\hbox{{${\buildrel
    \ldotp\over{\hbox to 6pt{\hrulefill}}}$}}\;}

\newcounter{savenumi}

\newtheorem{theoremfoo}{Theorem}[section] 

\newtheorem{lemmafoo}[theoremfoo]{Lemma}

\newtheorem{conjecturefoo}[theoremfoo]{Conjecture}

\newtheorem{conventionfoo}[theoremfoo]{Convention}

\newtheorem{porismfoo}[theoremfoo]{Porism}

\newtheorem{gamefoo}[theoremfoo]{Game}

\newtheorem{corollaryfoo}[theoremfoo]{Corollary}

\newtheorem{claimfoo}[theoremfoo]{Claim}

\newtheorem{openfoo}[theoremfoo]{Open Problem}
\newenvironment{open}{\pagebreak[1]\begin{openfoo}\rm}{\end{openfoo}}

\newtheorem{exercisefoo}{Exercise}

\newcommand{\fig}[1] 
{
 \begin{figure}
 \begin{center}
 \input{#1}
 \end{center}
 \end{figure}
}

\newtheorem{potanafoo}[theoremfoo]{Potential Analogue}

\newtheorem{notefoo}[theoremfoo]{Note}

\newtheorem{notabenefoo}[theoremfoo]{Nota Bene}

\newtheorem{nttn}[theoremfoo]{Notation}

\newtheorem{empttn}[theoremfoo]{Empirical Note}

\newtheorem{examfoo}[theoremfoo]{Example}
\newenvironment{example}{\pagebreak[1]\begin{examfoo}\rm}{\end{examfoo}}

\newtheorem{dfntn}[theoremfoo]{Def}
\newenvironment{definition}{\pagebreak[1]\begin{dfntn}\rm}{\end{dfntn}}

\newtheorem{propositionfoo}[theoremfoo]{Proposition}

\newcommand{\yyskip}{\penalty-50\vskip 5pt plus 3pt minus 2pt}
\newcommand{\blackslug}{\hbox{\hskip 1pt
        \vrule width 4pt height 8pt depth 1.5pt\hskip 1pt}}
\newcommand{\QED}{{\penalty10000\parindent 0pt\penalty10000
        \hskip 8 pt\nolinebreak\blackslug\hfill\lower 8.5pt\null}
        \par\yyskip\pagebreak[1]}

\newcommand{\BBB}{{\penalty10000\parindent 0pt\penalty10000
        \hskip 8 pt\nolinebreak\hbox{\ }\hfill\lower 8.5pt\null}
        \par\yyskip\pagebreak[1]}

\newtheorem{factfoo}[theoremfoo]{Fact}

\begin{document}

\author{Johnathan Cai,  Ryan Diehl, William Gasarch,  Ian Kim,\\
Rohan Sinha}
\title{Estimating the Number of Primes In Unusual Domains} 

\maketitle

\noindent

\begin{abstract}

The Prime Number Theorem states that the number of primes in 
$\{1,\ldots,x\}$, denoted $\pi(x)$, is approximately 
$\frac{x}{\ln(x)}$. In this paper, we investigate the distribution of primes for 
domains other than $\N$. 
First we look at 
$\Ad=\{ x \colon x\equiv 1 \pmod d\}$. 
We give a heuristic argument to form a conjecture on the number of \textit{congruence monoid primes} in
$\Ad$ that are $\le x$.
We then provide empirical evidence that indicates our conjecture is close but may
need some correction. 
Second, we do similar calculations for the Gaussian Integers. 
Third, we discuss the difficulty of these types of questions for quadratic extensions of $\Z$.
\end{abstract}


\section{Introduction}

In 1896, Jacques Hadamard and Charles-Jean de la Vallée Poussin both independently discovered that the 
number of primes $\le x$ was roughly $\frac{x}{\ln(x)}$ using complex analysis. 
This was built off of earlier work done by Pafnuty Chebyshev in the 1850s. Later proofs have been 
created for the Prime Number Theorem by Paul Erdős and Atle Selberg in 1948. 
Their proof only uses calculus; however, they are still difficult. 
In 1980 Newman~\cite{Newman-1980} (also see Zagier~\cite{Zagier-1997}) gave a short proof that only
used a little complex analysis. 

In Section~\ref{se:def}, we define the {\it Congruence Monoid} and 
discuss an analog of the Prime Number Theorem for them. 
In Section~\ref{se:id}, we define the {\it Integral Domain} and related notions. 
In Section~\ref{se:gauss}, we define the Gaussian Integers and discuss
analogs of the prime number theorem in them.
In Section~\ref{se:other}, we discuss other integral domains. 
In Section~\ref{se:openprob}, we recap the open problems encountered.

\section{Congruence Monoids} \label{se:def} 

\subsection{Definitions For Congruence Monoids} 

\begin{definition}[Prime] 
Let \( \D \subseteq\ \N \). A \textit{prime} in $\D$ is a number $p \in \D$ such that 
$p > 1$  and the only positive divisors of $p$ that lie in $\D$ are $1$ and $p$ itself. 
\end{definition}

\begin{definition}[Congruence Monoid Prime]
Let \( d \in \N \) with \( d \geq 2 \). Define the set
\[
\Ad = \{ n \in \N \colon  n \equiv 1 \pmod{d} \}.
\]
An element \( p \in \Ad \) is called a \textit{congruence monoid prime} if $p\ne 1$ and, whenever \( p = ab \) for some \( a, b \in \Ad \), then either \( a = 1 \) or \( b = 1 \). Factorizations involving elements outside \( \Ad \) are not considered. 
\end{definition}

\noindent
{\bf Examples} 
Lets look at $\Af$. 
\begin{enumerate}
\item
We write down the first few elements:

$$1,5,9,13,17,21,25,29,33,37,41,45$$

\item
All of the numbers that are prime in $\N$ are prime in $\Af$.
So thats $5,13,17,29,37$.

\item
What about the numbers that are not prime in $\N$?

9 is prime. Even though $9=3\times 3$ note that $3\notin \Af$. 

21 is prime. Even though $21=3\times 7$, note that $e\notin\Af$. 

25 is not prime since $25=5\times 5$ and $5\in\Af$. 

33 is not primes since $33=3\times 11$ and $3\notin\Af$. 

45 is not prime since $45=5\times 9$ and $5,9\in\Af$. 

\item
Note that there are numbers that are not primes in $\N$ but are primes in 
$\Af$. 

\end{enumerate}

We want to define an analog of $\pi(x)$ for $\Ad$. We first carefully define $\pi(x)$.

\begin{definition}
$\pi(x)$ is the number of elements of $\{1,\ldots,x\}$ that are prime, i.e.,

$$\pi(x) = \#\{ x \in \{1,\ldots,x\}  \colon p \text{ is prime} \}.$$

\end{definition}

Note that the domain of interest is $\{1,\ldots,x\}$. 
We define $\pi_d(x)$, the analog of $\pi(x)$ for $\Ad$, noting that
the domain of interest is $\Ad \cap \{1,\ldots,x\}$. 

\begin{definition}[Congruence Monoid Prime Count]
$\pi_d(x)$ is the number of elements of $\Ad \cap \{1,\ldots,x\}$ that are prime in $\Ad$, i.e.,

$$\pi_d(x) = \#\{ x \in \Ad \cap \{1,\ldots,x\}  \colon p \text{ is prime in $\Ad$} \}.$$
\end{definition}

\subsection{Congruence Monoid Prime Estimation} \label{se:mono} 

How does $\pi_d(x)$ compare to $\pi(x)$? 

\begin{itemize}
\item
The domain for $\pi_d(x)$ is smaller than that for $\pi(x)$. 
This suggests that $\pi_d(x)\le \pi(x)$.
\item
There are primes in $\Ad$ that are not primes in $\N$. Hence
this suggests that $\pi(x)\le \pi_d(x)$.
\end{itemize}

\noindent Balancing these effects, we propose the estimation
\[
\pi_d(x) \approx \frac{x}{d \ln (x)^{1/d}}.
\]

Our empirical results from computational simulations support this approximation, and its accuracy can be quantified by the normalized ratio

\[
R_d(x) = \frac{\pi_d(x)}{\,x / \bigl(d\ln (x)^{1/d}\bigr)}.
\]
Values $R_d(x) \approx 1$ indicate close agreement between the model and observed data.

The table below displays the accuracy of the estimation across several values of \( d \) for primes \(\le\) $x=10^4$. The mean absolute percentage deviation (MAPE) measures the accuracy of a forecasting model, which in this case is our estimation. The corresponding graphs provide a visual comparison of \( \pi_d(x) \) and the estimation over the full range of \( x \) for a few values of \(d\) from the table. Complete data can be viewed \href{https://github.com/newrohansinha/Primes/blob/main/n%20%3D%2010%5E4.pdf}{here}.

\begin{table}[H]
\centering
\begin{tabular}{|c|r|r|r|r|r|r|}
\hline
$d$ & Largest Prime & Actual Count & Estimate & $R_d$ & $|R_d - 1|$ & MAPE (\%) \\
\hline
3  & 10000  & 1380 & 1590.21 & 0.86781 & 0.13219 & 9.05 \\
5  & 9996   & 1210 & 1282.34 & 0.94358 & 0.05642 & 3.81 \\
7  & 9997   & 1009 & 1039.97 & 0.97022 & 0.02978 & 2.45 \\
9  & 10000  & 851  & 868.19  & 0.98020 & 0.01980 & 2.28 \\
11 & 10000  & 745  & 742.93  & 1.00279 & 0.00279 & 2.88 \\
13 & 9998   & 653  & 648.33  & 1.00720 & 0.00720 & 3.10 \\
21 & 9997   & 438  & 428.29  & 1.02268 & 0.02268 & 3.88 \\
50 & 9951   & 196  & 190.38  & 1.02953 & 0.02953 & 3.03 \\
\hline
\end{tabular}
\caption{Comparison of Actual and Estimated D\textsubscript{d}-Prime Counts up to $10^4$}
\end{table}


\begin{figure}[H]
\centering
\includegraphics[width=0.8\textwidth]{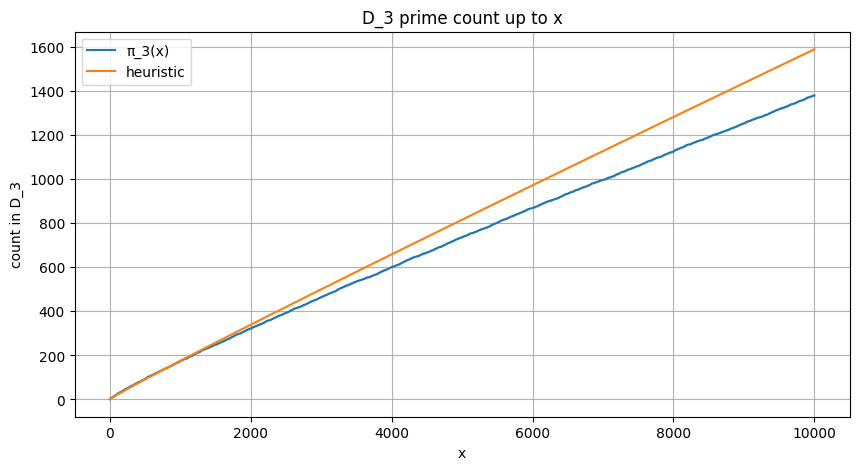}
\caption{Actual prime count \(\pi_{3}(x)\) (blue) versus the estimate (orange) for \(d = 3\).}
\label{fig:d3}
\end{figure}

The estimation very minorly underestimates until $x\approx800$, then increasingly starts to overestimate.

\vspace{1em}

\begin{figure}[H]
\centering
\includegraphics[width=0.8\textwidth]{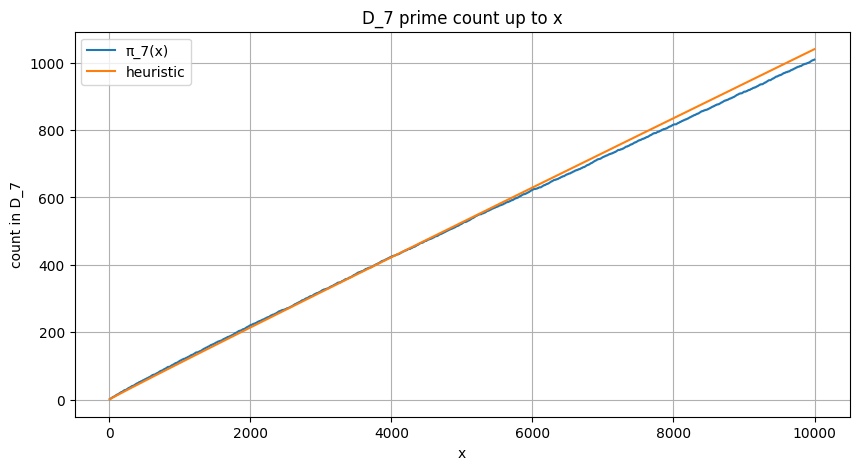}
\caption{Actual prime count \(\pi_{7}(x)\) (blue) versus the estimate (orange) for \(d = 7\).}
\label{fig:d7}
\end{figure}

Similar to figure 1, the estimate very minorly underestimates until $x\approx4100$ and then increasingly starts to overestimate.

\vspace{1em}

\begin{figure}[H]
\centering
\includegraphics[width=0.8\textwidth]{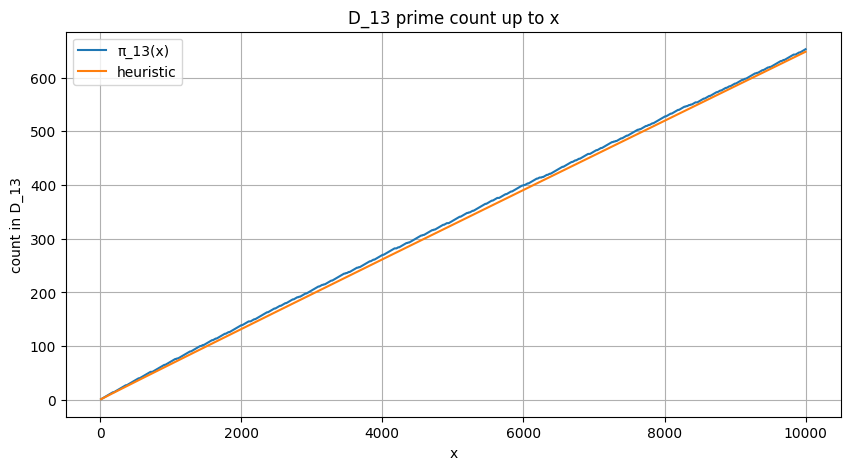}
\caption{Actual prime count \(\pi_{13}(x)\) (blue) versus the estimate (orange) for \(d = 13\).}
\label{fig:d13}
\end{figure}

Same phenomenon described above except at $ x\approx13800$ (not shown in graph).

\vspace{1em}

\begin{figure}[H]
\centering
\includegraphics[width=0.8\textwidth]{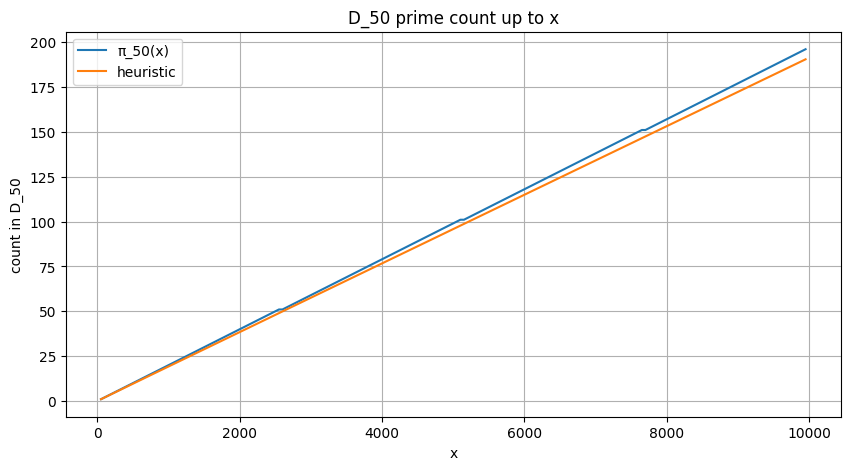}
\caption{Actual prime count \(\pi_{50}(x)\) (blue) versus the estimate (orange) for \(d = 50\).}
\label{fig:d50}
\end{figure}

Same phenomenon described above except at $ x\approx330,000$ (not shown in graph).

Unlike the Prime Number Theorem, where the estimation becomes more accurate as $x$ grows (in fact, becomes perfect as $x$ approaches infinity), the accuracy of our estimation is not as straightforward. The estimation becomes increasingly accurate up to a certain value of \(x\) (almost exact) and then diverges from the actual count.

Our conjecture seems to be, as they used to say, {\it close but no cigar} 
(note that none of the authors smoke and three of them are high school students). 

\begin{open}~
\begin{enumerate}
\item
Modify our conjecture so that it better fits the data.
\item
Prove that modification. Since the Prime Number Theorem was hard
to prove this might be hard as well; however, the math needed for the
Prime Number Theorem is known and may be useful. 
\end{enumerate}
\end{open}

\section{Integral Domains}\label{se:id}

When studying primes in integral domains careful definitions of
{\it unit}, {\it prime}, and  {\it irreducible} are needed.

\begin{definition}
An integral domain is a triple $\D=(D,+,\times)$ such that 
(1) $D$ is a set,
(2) $+$ and $\times$ are maps from $D\times D$ into $D$, and
(3) $(D,+,\times)$ satisfies the following properties:
\begin{enumerate}
\item
$D$ is closed under $+$ and $\times$.
\item
$+$ and $\times$ are both commutative and associative.
\item
$\times$ is distributive over $+$. Hence, for all $a,b,c\in D$,
$a\times (b+c) = a\times b + a\times c$. 
\item
There exists an element $0\in D$ such that, for all $a\in D$, $a+0=0+a=a$. 
\item 
For every $a\in D$ there exists $b\in D$ such that $a+b=0$. We usually
denote $b$ by $-a$. 
\item
There exists an element $1\in D$ such that, for all $a\in D$, $a\times 1=a$. 
\item 
For all $a,b\in D$, if $ab=0$ then either $a=0$ or $b=0$. 
\end{enumerate}
\end{definition}

\begin{definition}\label{de:id}
Let $\D$ be an integral domain.
\begin{enumerate}
\item 
An element $u\in \D$ is a \textit{unit} there exists  $v\in\D$ such that $uv=1$. 
The units of $\Z$ are $\{-1,1\}$.
\item 
An element $p\in\D$ is  \textit{prime (in $\D$)} if (a) $p$ is not a unit, and
(b) if $p$ divides $ab$ then either $p$ divides $a$ or $p$ divides $b$.
The primes of $\Z$ are the usual primes and their negations. For example, 3 and $-7$ are
both primes in $\Z$. 
\item 
An element $r \in \D$ is \textit{irreducible (in $\D$)} if (a) $r$ is not a unit, and 
(b) if $r= ab$, then either $a$ or $b$ is a unit.
In $\Z$ irreducibles and primes are the same. This equivalence is  false in some other integral
domains. 
\end{enumerate}
\end{definition}

\begin{definition}
Let $\alpha\notin\Q$. Then $\Z[\alpha]$ is the set $\{ a+b\alpha \colon a,b\in \Z\}$
\end{definition}

\section{Gaussian Prime Estimation}\label{se:gauss}

\begin{definition}
The {\it Gaussian Integers} is the set $\Z[i]$ where $i=\sqrt{-1}$.
\end{definition}

The Gaussian Integers are an integral domain. Hence the definitions of units, primes, and
irreducible from Definition~\ref{de:id} apply to them. 
In the Gaussian integers
(a) the units are $\{1,-1,i,-i\}$, and
(b) primes and irreducibles are the same. 

We want to study an analog of the Prime Number Theorem for the Gaussian Integers. Recall that the usual Prime Number Theorem is about the number of primes in 
$\{1,\ldots,x\}$. Since $\N$ is not an integral domain, but $\Z$ is, lets rewrite that as the prime number
theorem taking into account that the domain is $\Z$:

{\it The number of primes in $\{ y\in \Z \colon 0\le |y|\le x \}$ is approximately $\frac{x}{\ln(x)}$.} 
We take $0\le |y|$ since in the usual Prime Number Theorem we do not count both a prime and its
negation. For example, we don't count both 7 and $-7$ as primes.

For the Gaussian Integers we need (1) a notion of size analogous to absolute value for $\Z$,
(2) a way to not count $p$, $-p$, $ip$, and $-ip$. 

\begin{definition}~
\begin{enumerate}
\item 
Let $a,b\in\Z$. Then the \textit{norm} of $a+bi$ is $a^2 + b^2$.
This is denoted by $N(a+bi)$. 
This will not be our notion of size; however, $\sqrt{N(a+bi)}$ will be. 
\item 
Let $r\in\N$. The \textit{norm circle} of radius $r$ in the set

$$\{ a+ bi \colon \sqrt{a^2+b^2}\le r \}.$$
\end{enumerate}
\end{definition}

Using the norm circle provides a clear and finite boundary, making it possible to study the distribution of 
Gaussian primes up to a specific size, analogous to counting up to $x$ in the Prime Number Theorem.
There is one more issue: we only look at $a,b\ge 0$ since that avoids the problem of counting
a prime four times because of units. 

\begin{definition}
$\pi_G(r)$ is the number of primes in $a+bi\in \Z[i]$ such that
(a) $a,b\ge 0$, and (b) $a+bi$ is in the norm circle of radius $r$. 
\end{definition}

We empirically found the estimate $\pi_G(r)\approx \frac{r^2}{2\ln r}$

\begin{figure}[H]
\centering
\includegraphics[width=0.7\textwidth]{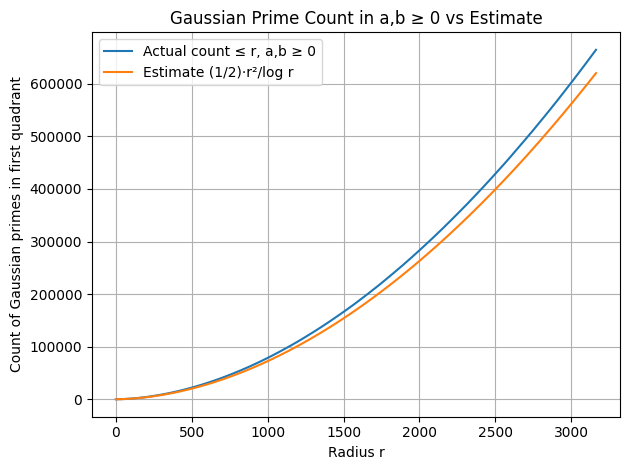}
\caption{Actual Gaussian prime count for \(
0 < r^{2} \le 10^7 
\) (blue) versus the the estimate (orange).}
\label{fig:gauss10m}
\end{figure}

\begin{table}[H]
\centering
\begin{tabular}{|c|r|}
\hline
Radius (r) & MAPE (\%) \\
\hline
$10^3$ & 18.630 \\
$10^4$ & 14.259 \\
$10^5$ & 10.907 \\
$10^6$ &  8.695 \\
$10^7$ &  7.220 \\

\hline
\end{tabular}
\caption{Mean absolute percentage deviation (MAPE) of the estimate \(\pi_G(r)\) at various radii.}
\label{tab:gauss-mape}
\end{table}

\begin{open}
Prove or disprove the estimate $\pi_G(r) \sim \frac{r^2}{2\ln r}$. 
If the estimate is true then find the error term. 
\end{open}

\section{Other Integral Domains}\label{se:other} 

We list integral domains and discuss if analogs of the Prime Number Theorem can
be asked and what issues might arise.
Much of the information in this section is from Weintraub~\cite{Weintraub-2008}. 

\noindent
1) $\Z[\sqrt{-d}]$ where $d\in\N$ is squarefree. The standard norm is 
$N(a+b\sqrt{-d})=a^2+db^2$. 
There are two possible questions to ask:
\begin{itemize}
\item 
Approximate the
number of primes of the form  $a+b\sqrt{-d}$ where $a,b\ge 0$ and
$a^2+db^2\le r$, or 
\item 
Approximate the
number of primes of the form  $a+b\sqrt{-d}$ where $a,b\ge 0$ and
$a^2+b^2\le r$. 
\end{itemize}

\bigskip

\noindent
2) $\Z[\sqrt{d}]$ where $d\in\N$ is squarefree. Asking about primes in
this integral domain is problematic for several reasons.

\noindent
2a) For some values of $d$ this is the wrong question.

In Section~\ref{se:gauss} we looked at $\Z[i]$. Why not $\Z[\frac{i}{2}]$?
Lets start with $\Q(i)$. We look for an integral domain $\D$ such that
$\Z[i]\subseteq \D\subseteq \Q(i)$ and 
(roughly): $\Q$ is to $\Z$ as $\Q(i)$ is to $\D$. 

\begin{definition}
Let $\D_1$ and $\D_2$ be integral domains such that $\D_1\subseteq \D_2$. 
\begin{enumerate}
\item
$x\in \D_2$ is {\it integral over $\D_1$} if $x$ is the root of a monic polynomial with
coefficients in $\D_1$. 
\item
The set of elements of $\D_2$ that are integral over $\D_1$ is the
{\it integral closure of $\D_1$ in $\D_2$}.
\end{enumerate}
\end{definition}

\begin{example}~
\begin{enumerate}
\item
$\Z$ is the integral closure of $\Z$ in $\Q$. 
\item
$\Z[i]$ is the integral closure of $\Z$ in $\Q[i]$. 
\item
$\Z[\sqrt{5}]$ is {\it not} the integral closure of $\Z$ in $\Q[\sqrt{5}]$:
$\frac{1+\sqrt{5}}{2}$ is integral over $\D$---it satisfies $x^2-x-1=0$. 
\item
The integral closure of $\Z$ in $\Q[\sqrt{5}]$ is $\Z[\frac{1+\sqrt{5}}{2}]$. 
\item
Let $d$ be a square free integer. 
\begin{enumerate}
\item
If $d\equiv 2,3 \pmod 4$ then the integral closure of $\Z$ in $\Q(\sqrt{d})$ is $\Z[\sqrt{d}]$. 
\item 
If $d\equiv 1 \pmod 4$ then the integral closure of $\Z$ in $\Q(\sqrt{d})$ is 
$\Z[\frac{1+\sqrt{5}}{2}]$. 
\end{enumerate}
\end{enumerate}
\end{example}

\begin{definition}
Let $d$ be a squarefree integer.
We define $\O(\sqrt{d})$ as follows:
\begin{enumerate}
\item 
If $d\equiv 2,3 \pmod 4$ then $\O(\sqrt{d})=\Z[\sqrt{d}]$. 
\item
If $d\equiv 1 \pmod 4$ then $\O(\sqrt{d})=\Z[\frac{\sqrt{d}+1}{2}]$. 
\end{enumerate}
\end{definition}

Here is the right question to ask: 

{\it is there an analog of the Prime Number Theorem for $\O(d)$?} 

\bigskip

\noindent
2c) Number of Units.

\begin{enumerate}
\item
If $d\ge 1$ then $\O(\sqrt{d})$ has an infinite number of units. This makes it hard
to phrase an analog of the Prime Number Theorem. 
\item
$\O(\sqrt{-1})$ has four units ($\pm 1$, $\pm i$). $\O(\sqrt{-3})$ has six units
(the six roots of unity). For all squarefree naturals $d\ge 5$, $\O(\sqrt{-d})$ has
two roots of unit ($\pm 1$). Hence for these an analog of the Prime Number
Theorem may be possible. But see the next item.
\end{enumerate}

\bigskip

\noindent
2d) For some $\O(\sqrt{d})$ primes and irreducibles are not the same. This may cause problems.

\bigskip

\noindent
2e) For some $d$ $\O(\sqrt{d})$ is not a unique factorization domain. This may cause problems.

\bigskip

We are {\it not} saying that formulating an analog of the Prime Number Theorem
in $\O(\sqrt{d})$ is impossible; however, there are some difficulties to overcome. 

\bigskip

\noindent
3) What about adding cube-roots or higher fractional powers? What about adding more
irrationals? These get into issues far harder than those encountered
for quadratic extension. 

\section{Open Problems} \label{se:openprob}

We recap the open problems stated earlier. 

\begin{enumerate}
\item
Prove or disprove the suggested estimate for the number of primes in $A_d$. 
If proven then obtain an error term. 
If disproven the find the correct approximation. 

\item
Prove or disprove the suggested estimate for the number of primes in $\Z[i]$. 
If proven then obtain an error term. 
If disproven the find the correct approximation. 

\item
Formulate analogs of the Prime Number Theorem in $\Z[\sqrt{d}]$ for various
values of $d$. Get empirical evidence to formulate conjectures. 
\end{enumerate}


\end{document}